\documentclass[a4paper,11pt,twoside,reqno]{article}

\usepackage{amsmath,amsthm,amsfonts,latexsym,amscd,amssymb}
\usepackage[a4paper, total={6in,8in}]{geometry}
\numberwithin{equation}{section}
\input xypic
\def\qed{{\hbadness=10000\hfill\ \vbox{\hrule height.09ex
			\hbox{\vrule width.09ex height1.55ex depth.2ex \kern1.8ex
				\vrule width.09ex height1.55ex depth.2ex}\hrule height.09ex}\break
		\bigskip}}
\newtheorem{theorem}{Theorem}[section]
\newtheorem{lemma}{Lemma}[section]
\newtheorem{corollary}{Corollary}[section]
\newtheorem{proposition}{Proposition}[section]

\theoremstyle{definition}
\newtheorem{definition}{Definition}[section]

\theoremstyle{remark}
\newtheorem{remark}{Remark}[section]
   
\newcommand{\n}{\noindent}

\begin{document}
	
	\linespread{1}\title {\textbf{Almost Gradient Ricci Solitons on Static Spacetime}}
	
	\author{Akhilesh Yadav, Tarun Saxena\thanks{Corresponding author}\\{\n E-mail: tarunsaxena254@gmail.com\thanks{Corresponding author E-mail}}}
	\date{}
	
	\maketitle 
	
	\noindent\textbf{Abstract:} The aim of this paper is to study geometrical aspects of static spacetime admitting an almost gradient Ricci soliton. Among others, We first determine the conditions under which the base manifold of static spacetime possess an almost gradient Ricci soliton and we show that the almost gradient Ricci soliton become steady gradient Ricci soliton when static spacetime turns to a vacuum static spacetime. Next, we exhibit that an expanding almost gradient Ricci soliton on base manifold of non-compact and connected static spacetime satisfies shr$\ddot{o}$dinger's equation for a smooth function $f$. Also, we find the soliton constant under which the static perfect fluid spacetime with almost gradient Ricci soliton holds the null convergence condition and the strong energy condition. Further, we study the almost gradient Ricci soliton on base manifold of static perfect fluid spacetime with potential function as warping function and it is shown that the base manifold of a static perfect fluid spacetime with an almost gradient Ricci soliton is an Einstein manifold. Next, we obtain a necessary and sufficient condition on soliton constant to obey timelike convergence condition. Further, we obtain some results for Ricci symmetric and weakly Ricci symmetric base manifold of static perfect fluid spacetime admitting gradient Ricci soliton. Finally, we find the nature of almost gradient Ricci soliton on $4$-dimensional half conformally flat base manifold of static perfect fluid spacetime.

	\n\textbf {Mathematics Subject Classification (2020):} 53C50, 53Z05, 83C05, 83C56.

	\n\textbf{Key words:} Static spacetime, Almost Gradient Ricci soliton, Perfect fluid spacetime, Shr$\ddot{o}$dinger's equation, Ricci symmetric spacetime, Stiff matter.

	\section{Introduction}
	In recent years, research interest in Ricci solitons has increased significantly. The concept of Ricci solitons was introduced by Hamilton $[13]$ as a natural generalization of an Einstein metric. A Riemannian metric $g$ on a smooth manifold $M$ is said to be a Ricci soliton, if there exist a constant $\lambda$, and a smooth vector field $V$ on $M$ satisfying
	\begin{equation}
		(\,L_V g)\,(\,X_1, X_2)\, + 2S(\,X_1, X_2)\, + 2\lambda g(\,X_1, X_2)\, = 0,
	\end{equation}
	for all smooth vector fields $X_1, X_2$ on $M$, where $L_V g$ denotes the Lie derivative of $g$ along the vector field $V$ on $M$ and $S$ is the Ricci tensor of $g$. Here, the vector field $V$ and the constant $\lambda$ are known as potential vector field and the soliton constant, respectively. The Ricci soliton is called shrinking, steady, and expanding if $\lambda < 0$, $\lambda = 0$, and $\lambda > 0$, respectively. Obviously, Einstein metrics are the trivial Ricci solitons with potential vector field $V$ as homothetic vector field or zero. Ricci solitons are the fixed points of the Hamilton's Ricci flow $[13]$
	\begin{equation}
		\frac{\partial g(t)}{\partial t} = -2S(g(t)),
	\end{equation}
	viewed as a dynamical system on the space of Riemannian metrics modulo diffeomorphisms and scaling. Also, the Ricci solitons model the formation of singularities in the Ricci flow and correspond to self-similar solutions $[17]$.
	
	\n If the potential vector field $V$ in $(1.1)$ is considered as the gradient of some smooth function $f$ on $M$, then the Ricci soliton is called a gradient Ricci soliton, and $(1.1)$ can be written as
	\begin{equation}
		H^f(X_1, X_2) + S(X_1, X_2) + \lambda g(X_1, X_2) = 0,
	\end{equation}
	for all smooth vector fields $X_1, X_2$ on $M$, where $H^f$ is the Hessian of potential function $f$, $S$ is the Ricci tensor of $g$ and $\lambda$ is a constant on $M$. In $[19]$, S. Pigola generalized the notion of Ricci soliton to almost Ricci soliton by setting the soliton constant $\lambda$ to be a smooth function on $M$.
	
	\n An n-dimensional semi-Riemannian manifold together with a metric of signature $(1, (n - 1))$ is named as a Lorentzian manifold of dimension $n$ $[18]$. A general relativistic spacetime is regarded as a connected four-dimensional Lorentzian manifold equipped with Lorentzian metric of signature $(-,+,+,+)$. For many years, a lot of interest has been noticed in the study of physical properties of spacetimes with some geometrical structures. As static spacetimes are considered as generalizations of the static vacuum spaces and hence static spactimes are important topic of research for both physicists and mathematicians. In $[14]$, Kobayashi and Obata studied static spacetime on a Lorentzian manifold obeying Einstein field equation with perfect fluid as a perfect matter, and discussed conformally flatness of such spacetime. In $[7]$, authors obtained conditions for the Riemannian factor and the warping function of a standard static spacetime. In $[10]$, authors studied static perfect fluid spacetimes on GRW spacetimes, and under certain restrictions they proved such a spacetime is of Petrov type I, D or O. Many authors discussed the characterization of static perfect fluid spacetimes on contact metric manifold $[5]$, compact manifolds $[6]$, half conformally flat $[15]$, and paracontact metric manifolds $[20]$. In $[25]$, authors studied the idea of Ricci solitons of spherically symmetric static spacetimes, and found that special classes of such spacetime metrics admit shrinking, expanding or steady Ricci solitons. Also, in $[23]$, the authors studied standard static spacetime in terms of almost Ricci-Yamabe soliton with conformal vector field and discuss the harmonic aspects on standard static spacetime.
	
	\n After developing spaces of constant curvature, mathematicians classified locally symmetric Riemannian manifolds. The properties of locally symmetric manifolds have been studied by several authors in perfect fluid spacetime such as in $[16]$ Mallick and De studied some conditions for the existence of perfect fluid semi-Riemannian symmetric spacetime. Also, authors discussed pseudo symmetric manifolds $[4]$,  weakly symmetric manifolds $[26]$. This scenario inspired several mathematician and physicists to develop study of different Ricci solitons on spacetimes such as Ricci soliton ($[24]$, $[27]$), $\eta$-Ricci-Bourguignon soliton $[11]$, Ricci Yamabe solitons $[22]$, and conformal Ricci soliton $[21]$.
	
	\n Inspired by the above studies, here we consider static spacetime whose metric is an almost gradient Ricci soliton. The present paper is organized as follows: In sections 1 and 2, we write some basic information about Ricci solitons, static spacetimes, weakly Ricci symmetric spacetime, perfect fluid spacetimes, which are required for this paper. In section 3, we first determine the conditions under which the base manifold of static spacetime possess an almost gradient Ricci soliton and we show that the almost gradient Ricci soliton become steady gradient Ricci soliton when static spacetime turns to a vacuum static spacetime. Next, we exhibit that an expanding almost gradient Ricci soliton on base manifold of non-compact and connected static spacetime satisfies shr$\ddot{o}$dinger's equation for a smooth function $f$. In section 4, we find the soliton constant under which the static perfect fluid spacetime with almost gradient Ricci soliton holds the null convergence condition and the strong energy condition. In section 5, we study the almost gradient Ricci soliton on base manifold of static perfect fluid spacetime with potential function as warping function and it is shown that the base manifold of a static perfect fluid spacetime with an almost gradient Ricci soliton is an Einstein manifold. Next, we obtain a necessary and sufficient condition on soliton constant to obey timelike convergence condition. Further, we obtain some results for Ricci symmetric and weakly Ricci symmetric base manifold of static perfect fluid spacetime admitting gradient Ricci soliton. Finally, we find the nature of almost gradient Ricci soliton on $4$-dimensional half conformally flat base manifold of static perfect fluid spacetime.
	
	\section{Preliminaries}
	An $(n + 1)$-dimensional Lorentzian warped product manifold $\tilde{M}^{(n + 1)} = M^n \times_f \Re$, is called a static spacetime $[14]$ with the metric
	\begin{equation}
		\tilde{g} = g - f^2 dt^2,
	\end{equation}
where $M$ is an $n$-dimensional Riemannian manifold with metric $g$ and $f(>0)$ is a smooth function on $M$. We consider Einstein field equation on $(\tilde{M}, \tilde{g})$ with perfect fluid as a matter field, given by
	\begin{equation}
		\tilde{S} - \frac{\tilde{r}}{2}\tilde{g} = -\rho g -\sigma f^2 dt^2,
	\end{equation}
where $\tilde{S}$ and $\tilde{r}$ represent the Ricci tensor and the scalar curvature of $\tilde{g}$, respectively. Also, $\rho$ and $\sigma$ are smooth functions on $\tilde{M}$, which are called isotropic pressure and energy density, respectively. From $(2.1)$ and $(2.2)$,  the static perfect fluid (or simply, SPF) spacetime is equivalent to
	\begin{equation}
		f\left(S - \frac{r}{n}g\right) = H^f - \frac{\triangle f}{n}g,
	\end{equation}
and
  \begin{equation}
  	\triangle f = \left(\frac{n - 2}{2(n - 1)}r + \frac{n}{n - 1}\rho\right)f,
  \end{equation}
where $S$ and $r$ denotes the Ricci tensor and the scalar curvature on $M^n$, respectively and $\triangle f$ is the Laplacian of smooth function $f$, and hence, the static spacetime $(\tilde{M}, \tilde{g})$ is considered as static perfect fluid (SPF) spacetime. In particular, if $\rho = r = 0$ in $(2.3)$ and $(2.4)$, then the static perfect fluid equation reduces in static vacuum Einstein equations given by
\begin{equation}
	fS = H^f, \quad  \triangle f = 0.
\end{equation}
Also, if we put $\rho = -\frac{r}{2}$ in $(2.4)$ and using it in $(2.3)$, we have an equation, named Fischer-Marsden equation given by
\begin{equation}
	H^f - fS - \triangle f g = 0.
\end{equation}
Let $\tilde{M} = M \times_f \Re$ with the metric $\tilde{g} = g - f^2 dt^2$, be a static spacetime and $\nabla$ be the Levi-Civita connection on Riemannian manifold $M$. Then the Levi-Civita connection $\tilde{\nabla}$ on $\tilde{M}$ is given by $[23]$
\begin{equation}
	\tilde{\nabla}_{\partial_t}\partial_t = fgradf, \quad \tilde{\nabla}_{\partial_t}X_1 = \tilde{\nabla}_{X_1}\partial_t = X_1(\ln f)\partial_t, \quad \tilde{\nabla}_{X_1}X_2 = \nabla{X_1}X_2,
\end{equation}
for all smooth vector fields $X_1, X_2$ on $M$. Let $R$ be the Riemannian curvature tensor, $S$ be the Ricci curvature tensor and $r$ be the scalar curvature of $M$. Then Riemannian curvature tensor $\tilde{R}$, Ricci curvature tensor $\tilde{S}$ and scalar curvature $\tilde{r}$ of $(\tilde{M}, \tilde{g})$ are given by $[23]$
\begin{equation}
	\tilde{R}(X_1, \partial_t)\partial_t = -f\nabla_{X_1}gradf,
\end{equation}
\begin{equation}
	\tilde{R}(\partial_t, \partial_t)\partial_t = \tilde{R}(\partial_t, \partial_t)X_1 = \tilde{R}(X_1, X_2)\partial_t = 0,
\end{equation}
\begin{equation}
	\tilde{R}(\partial_t, X_1)X_2 = \frac{1}{f}H^f(X_1, X_2)\partial_t,
\end{equation}
\begin{equation}
	\tilde{R}(X_1, X_2)X_3 = R(X_1, X_2)X_3,
\end{equation}
\begin{equation}
	\tilde{S}(\partial_t, \partial_t) = f\triangle f,
\end{equation}
\begin{equation}
	\tilde{S}(X_1, \partial_t) = 0,
\end{equation}
\begin{equation}
	\tilde{S}(X_1, X_2) = S(X_1, X_2) - \frac{1}{f}H^f(X_1, X_2),
\end{equation}
\begin{equation}
	\tilde{r} = r -2\frac{1}{f}\triangle f,
\end{equation}
where $X_1, X_2, X_3$ are smooth vector fields on $M$, $H^f(X_1, X_2)$ is the Hessian of $f$ on $M$ and $\triangle f$ denotes the Laplacian of $f$ on $M$. 

\n The notion of weakly Ricci symmetric manifold was first used by Tamassy and Binh $[26]$ in 1993. A non-flat Riemannian manifold $(M, g)$ of dimension $n$ $(n > 2)$ is said to be weakly Ricci symmetric manifold, simply $(wRs)_n$- manifold if Ricci tensor $S$ is not identically zero and satisfying 
	\begin{equation}
		(\nabla_{X_1} S)(X_2, X_3) = \alpha(X_1)S(X_2, X_3) + \beta(X_2)S(X_3, X_1) + \gamma(X_3)S(X_1, X_2),
	\end{equation}
	for all smooth vector fields $X_1, X_2, X_3$ on $M$, where $\alpha$, $\beta$ and $\gamma$ are three non-zero 1-forms and $\nabla$ denotes the operator of covariant differentiation with respect to the metric $g$.
	On a $(wRs)_n$- manifold there is a 1-form $\delta$ defined as
	\begin{equation}\nonumber
		\delta(X) = \beta(X) - \gamma(X).
	\end{equation}
	\begin{lemma}$[8]$
		If $\delta \neq 0$ on a $(wRs)_n$- manifold $(M, g)$, then the Ricci tensor $S$ is defined by
		\begin{equation}
			S(X_1, X_2) = -r\mu(X_1)\mu(X_2),
		\end{equation}
		for all smooth vector fields $X_1, X_2$ on $M$, where $r$ is scalar curvature on $M$ and $\mu$ is a non-zero 1-form given by
		\begin{equation}
			\mu(X) = g(X, \rho),
		\end{equation}
		for all smooth vector field $X$ on $M$ and $\rho$ is called basic vector field.
	\end{lemma}
	
	\n A four-dimensional semi-Riemannian manifold $M$ with Lorentzian metric $g$ is called a weakly Ricci symmetric spacetime or simply $(wRs)_4$ spacetime, if its Ricci tensor $S$ satisfies $(2.16)$.
	\begin{lemma}$[8]$
		In a $(wRs)_4$ spacetime, the scalar curvature $r$ is non-zero. If $r = 0$, then from $(2.2)$ the Ricci tensor $S = 0$ which is not possible by the definition of $(wRs)_4$ spacetime.
	\end{lemma}
	\begin{definition}$[18]$
		An $n$ dimensional Lorentzian manifold $M$ whose Ricci tensor $S$ satisfies the condition
		\begin{equation}\nonumber
			S(X_1, X_2) = \alpha g(X_1, X_2) + \beta \eta(X_1)\eta(X_2),
		\end{equation}
		for all smooth vector fields $X_1, X_2$ on $M$, is said to be a perfect fluid spacetime; where $\alpha$ and $\beta$ being the scalar fields and $\rho$ is unit timelike vector field, called velocity
		vector field defined by $g(X, \rho) = \eta(X)$.
	\end{definition}
	\begin{definition}
		A smooth vector field $\xi$ on a Riemannian manifold $M$ is called a torseforming vector field if
		\begin{equation}
			\nabla_{X}\xi = fX + \nu(X)\xi,
		\end{equation}
		for all smooth vector field $X$ on $M$, where $f$ is a smooth function and $\nu$ is a 1-form on $M$.
	\end{definition}
\begin{lemma}$[3]$
	If an $n$-dimensional manifold $(M, g)$ is admitting an almost gradient Ricci soliton then the following identities hold:
	\begin{enumerate}
		\item[(1)] $\triangle$f + r + n$\lambda$ = 0,
		\item[(2)] $\nabla$$\triangle$f + $\nabla$r + n$\nabla$$\lambda$ = 0,
		\item[(3)] $\nabla$$\triangle$f + Ric($\nabla$f) + $\frac{1}{2}$$\nabla$r + $\nabla$$\lambda$ = 0,
		\item[(4)] $R(X_1, X_2, X_3, \nabla f) = d(\lambda)(X_2)(X_1, X_3) - d(\lambda)(X_1)(X_2, X_3) + (\nabla_{X_2}S)(X_1, X_3) - (\nabla_{X_1}S)(X_2, X_3)$,
	\end{enumerate}
for all smooth vector fields $X_1, X_2, X_3$ on $M$.
\end{lemma}
	
	\n The Weyl conformal curvature tensor $\mathcal{W}$ on the manifold $(M, g)$ is defined as follows:
	\begin{align}
		\mathcal{W}(X_1, X_2)X_3 =& R(X_1, X_2)X_3 + \frac{r}{(n - 1)(n - 2)}[g(X_1, X_3)X_2 - g(X_2, X_3)X_1] \\ \nonumber
		& -\frac{1}{(n - 2)}[S(X_2, X_3)X_1 - S(X_1, X_3)X_2 + g(X_2, X_3)QX_1 - g(X_1, X_3)QX_2],
	\end{align}
and the cotton tensor $\mathcal{C}$ is given by
\begin{align}
	\mathcal{C}(X_1, X_2)X_3 =& (\nabla_{X_1}S)(X_2, X_3) - (\nabla_{X_1}S)(X_2, X_3) \\\nonumber
	& \frac{1}{2(n - 1)}[g(X_2, X_3)dr(X_1) - g(X_1, X_3)dr(X_2)],
\end{align}
for all smooth vector fields $X_1, X_2, X_3$ on $M$, where $Q$ is the Ricci operator of metric $g$ on $M$. Using $(2.11)$, $(2.12)$, and $Lemma(2.3)$, we can write the expression of Weyl conformal curvature tensor as follows:
\begin{lemma}$[3]$
	If an $n$-dimensional manifold $(M, g)$ is admitting an almost gradient Ricci soliton then
	\begin{align}
		\mathcal{W}(X_1, X_2, X_3, \nabla f) =& \frac{r}{(n - 1)(n - 2)}[g(X_1, X_3)g(X_2, \nabla f) - g(X_2, X_3)g(X_1, \nabla f)] \\ \nonumber
		& \frac{1}{(n - 1)(n - 1)}\{g(X_2, X_3)S(X_1, \nabla f) - g(X_1, X_3)S(X_2, \nabla f)\} \\ \nonumber
		& \frac{1}{(n - 2)}\{S(X_2, X_3)g(X_1, \nabla f) - S(X_1, X_3)g(X_2, \nabla f)\} - \mathcal{C}(X_1, X_2, X_3),
	\end{align}
for all smooth vector fields $X_1, X_2, X_3$ on $M$.
\end{lemma}

\n We now quote some necessary propositions which will be further used in our theory.
\begin{proposition}$[1]$
	Let $M = (a, b)_f \times H$ be a standard static spacetime with
	$(H, h)$ a Ricci flat Riemannian manifold dim $H \geq 2$. Then $h(w, w)\triangle f - Hess(f )(w, w) \geq 0$ for all $w \in TH$ if and only if $M$ satisfies the null convergence condition.
\end{proposition}
\begin{proposition}$[1]$
	Let $M=(a,b)_f \times H$ be a standard static spacetime with $(H, h)$ Ricci flat and dim $H \geq 2$. Then $M$ satisfies the null convergence condition if and only if it satisfies the strong energy condition. 
\end{proposition}
	\section{Almost Gradient Ricci Solitons on Static spacetime}
	
	In this section, we study an almost gradient Ricci soliton structure on static spacetime $(\tilde{M}^{n + 1}, \tilde{g})$ with potential function as $\tilde{\phi}$ and soliton constant $\tilde{\lambda}$, defined as
	\begin{equation}
		H^{\tilde{\phi}}(\tilde{X_1}, \tilde{X_2}) + \tilde{S}(\tilde{X_1}, \tilde{X_2}) + \tilde{\lambda}\tilde{g}(\tilde{X_1}, \tilde{X_2}) = 0,
	\end{equation}
	for all smooth vector fields $\tilde{X_1}, \tilde{X_2}$ on $\tilde{M}$, where $H^{\tilde{\phi}}$ denotes the Hessian of smooth function $\tilde{\phi}$ on $\tilde{M}$, $\tilde{S}$ is the Ricci tensor of $\tilde{g}$ and $\tilde{\lambda}$ is a smooth function on $\tilde{M}$.
	\begin{theorem}
		Let $\tilde{M} =  M \times_f \Re$ be a static spacetime with metric $\tilde{g} = g - f^2 dt^2$. If $(\tilde{g}, \tilde{\lambda}, \tilde{\phi})$ is an almost gradient Ricci soliton on $\tilde{M}$ and $H^f = 0$, then $(g, \lambda, \phi)$ is an almost gradient Ricci soliton on $M$.
	\end{theorem}
\begin{proof}
	Let $(\tilde{M}^{n + 1}, \tilde{g})$ be a static spacetime, and admits an almost gradient Ricci soliton. Then computing $(3.1)$ at $(X_1, X_2)$ on $M$, we get
	\begin{equation}
		H^{\tilde{\phi}}(X_1, X_2) + \tilde{S}(X_1, X_2) + \tilde{\lambda}\tilde{g}(X_1, X_2) = 0,
	\end{equation}
for all smooth vector fields $X_1, X_2$ on $M$. It is known that$[4]$, $H^{\tilde{\phi}}(X_1, X_2) = H^\phi(X_1, X_2)$, where $\phi = \tilde{\phi}\vert_M$, and using $(2.1)$ and $(2.14)$, we obtain
\begin{equation}
	H^\phi(X_1, X_2) + S(X_1, X_2) + \tilde{\lambda}g(X_1, X_2) - \frac{1}{f}H^f(X_1, X_2) = 0,
\end{equation}
for all smooth vector fields $X_1, X_2$ on $M$. Suppose that $H^f = 0$, then $(3.3)$ provides
\begin{equation}
	H^\phi(X_1, X_2) + S(X_1, X_2) + \tilde{\lambda}g(X_1, X_2) = 0,
\end{equation}
for all smooth vector fields $X_1, X_2$ on $M$. Thus, $(g, \lambda, \phi)$ is an almost gradient Ricci soliton with $\tilde{\lambda} = \lambda$ on $M$. This completes the proof.
\end{proof}
\begin{theorem}
	Let $\tilde{M} =  M \times_f \Re$ be a static spacetime with metric $\tilde{g} = g - f^2 dt^2$. If $(\tilde{g}, \tilde{\lambda}, \tilde{\phi})$ is an almost gradient Ricci soliton on $\tilde{M}$, then $(g, \phi, \lambda, \mu)$ is an almost gradient $\eta$-Ricci soliton on $M$.
\end{theorem}
\begin{proof}
	By using the fact that $[2]$
	\begin{equation}\nonumber
		\frac{1}{f}H^f = H^{\ln f} + \frac{1}{f^2}df \otimes df,
	\end{equation}
with $(3.3)$ we get
\begin{equation}
	H^\phi(X_1, X_2) + S(X_1, X_2) + \tilde{\lambda}g(X_1, X_2) - H^{\ln f} + \frac{1}{f^2}df \otimes df,
\end{equation}
for all smooth vector fields $X_1, X_2$ on $M$. Let us define $\ln f = l$ and using $dl = \eta$ in $(3.5)$, we obtain
\begin{equation}
	H^\phi(X_1, X_2) + S(X_1, X_2) + \tilde{\lambda}g(X_1, X_2) - H^l(X_1, X_2) - \eta(X_1)\eta(X_2) = 0,
\end{equation}
which implies
\begin{equation}
	H^\psi(X_1, X_2) + S(X_1, X_2) + \tilde{\lambda}g(X_1, X_2) - \eta(X_1)\eta(X_2) = 0,
\end{equation}
for all smooth vector fields $X_1, X_2$ on $M$, where $\psi = \phi - l$. Thus, $(g, \phi, \lambda, \mu)$ is an almost gradient $\eta$-Ricci soliton with $\tilde{\lambda} = \lambda$ and $\mu = -1$ on $M$. This finishes the proof.
\end{proof}
\begin{theorem}
	Let $\tilde{M} =  M \times_f \Re$ be a vacuum static spacetime with metric $\tilde{g} = g - f^2 dt^2$. If $(\tilde{g}, \tilde{\lambda}, \tilde{\phi})$ is an almost gradient Ricci soliton on $\tilde{M}$ and $H^{\tilde{\phi}}(\partial_t, \partial_t) = 0$ then the almost gradient Ricci soliton is steady.
\end{theorem}
\begin{proof}
	Let $(\tilde{M}^{n + 1}, \tilde{g})$ be a static spacetime, and admits an almost gradient Ricci soliton. Then computing $(3.1)$ at $(\partial_t, \partial_t)$ on $M$, we get
	\begin{equation}
		H^{\tilde{\phi}}(\partial_t, \partial_t) + \tilde{S}(\partial_t, \partial_t) + \tilde{\lambda}\tilde{g}(\partial_t, \partial_t) = 0,
	\end{equation}
Using $(2.1)$ and $(2.12)$ in $(3.8)$, we obtain
\begin{equation}
	H^{\tilde{\phi}}(\partial_t, \partial_t) + f\triangle f - \tilde{\lambda}f^2 = 0.
\end{equation}
Suppose that $H^{\tilde{\phi}}(\partial_t, \partial_t) = 0$, then $(3.9)$ provides
\begin{equation}
	f\triangle f - \tilde{\lambda}f^2 = 0.
\end{equation}
Since $(\tilde{M}, \tilde{g})$ is a vacuum static spacetime, so using $(2.5)$ in $(3.10)$, we get
\begin{equation}
	\tilde{\lambda} = 0.
\end{equation}
Hence, the proof is completed.
\end{proof}
\begin{theorem}
	Let $\tilde{M} =  M \times_f \Re$ be a connected and non-compact static spacetime with metric $\tilde{g} = g - f^2 dt^2$ and $(\tilde{g}, \tilde{\lambda}, \tilde{\phi})$ be an almost gradient Ricci soliton on $\tilde{M}$. If $(g, \lambda, \phi)$ is an expanding almost gradient Ricci soliton on $M$, i.e. $\lambda > 0$, and $H^{\tilde{\phi}}(\partial_t, \partial_t) = 0$ then the smooth function $f : M \longrightarrow \mathfrak{R}$ satisfies the shr$\ddot{o}$dinger's equation.
\end{theorem}
\begin{proof}
	Since $(g, \lambda, \phi)$ is an almost gradient Ricci soliton on $M$, so Using $(3.10)$ with $theorem(3.1)$, we obtain
	\begin{equation}\nonumber
		\triangle f - \lambda f = 0,
	\end{equation}
which implies
\begin{equation}
	(\triangle - \lambda)f = 0.
\end{equation}
As $(g, \lambda, \phi)$ is an expanding almost gradient Ricci soliton on $M$ i.e., $\lambda > 0$, so by $[12]$ and $(3.12)$ the smooth function $f : M \longrightarrow \mathfrak{R}$ satisfies the shr$\ddot{o}$dinger's equation. This completes the proof.
\end{proof}

\section{Almost Gradient Ricci Solitons on Static Perfect Fluid spacetime}
In this section, we study an almost gradient Ricci soliton structure on static perfect fluid spacetime $(\tilde{M}^{n + 1}, \tilde{g})$ with potential function as $\tilde{\phi}$ and soliton constant $\tilde{\lambda}$.
\begin{theorem}
	Let $\tilde{M} =  M \times_f \Re$ be a static perfect fluid spacetime with metric $\tilde{g} = g - f^2 dt^2$. If $(\tilde{g}, \tilde{\lambda}, \tilde{\phi})$ is an almost gradient Ricci soliton on $\tilde{M}$, and $H^{\tilde{\phi}}(\partial_t, \partial_t) = 0$ then the soliton constant
	\begin{equation}\nonumber
		\tilde{\lambda} = n\rho + \left(\frac{n}{2} - 1\right)\tilde{r}.
	\end{equation}
\end{theorem}
\begin{proof}
	Since $\tilde{M} =  M \times_f \Re$ be a static perfect fluid spacetime with metric $\tilde{g} = g - f^2 dt^2$, so using $(2.4)$ in $(3.10)$, we obtain
	\begin{equation}
		\tilde{\lambda} = \frac{n - 2}{2(n - 1)}r + \frac{n}{n - 1}\rho.
	\end{equation}
Now, using $(2.15)$, $(3.10)$ with $(4.1)$, we find
\begin{equation}\nonumber
		\tilde{\lambda} = n\rho + \left(\frac{n}{2} - 1\right)\tilde{r}.
\end{equation}
This finishes the proof.
\end{proof}
\begin{theorem}
	Let $\tilde{M} =  M \times_f \Re$ be a static perfect fluid spacetime with $(M, g)$ a Ricci flat Riemannian manifold and dim$M \geq 2$, and $(\tilde{g}, \tilde{\lambda}, \tilde{\phi})$ is an almost gradient Ricci soliton on $\tilde{M}$ with $H^\phi = 0$, where $\phi = \tilde{\phi}\vert_M$ then $\tilde{M}$ satisfies the null convergence condition if 
	\begin{equation}\nonumber
		\tilde{\lambda} \leq -\frac{r}{n - 1}
	\end{equation}
\end{theorem}
\begin{proof}
	Let $(\tilde{g}, \tilde{\lambda}, \tilde{\phi})$ is an almost gradient Ricci soliton on static perfect fluid spacetime $\tilde{M} =  M \times_f \Re$. Then using $(2.3)$, we obtain
	\begin{equation}
		g(X, X)\triangle f - H^f(X, X) = (n - 1)H^f(X, X) - fnS(X, X) + rfg(X, X),
	\end{equation}
for all smooth vector fields $X$ on $M$. Since $(M, g)$ is a Ricci flat Riemannian manifold, then using $(3.3)$ in $(4.2)$ and $S(X, X) = 0$, we obtain
\begin{equation}
	g(X, X)\triangle f - H^f(X, X) = f(n - 1)H^\phi(X, X) + \{\tilde{\lambda}(n - 1) + r\}fg(X, X),
\end{equation}
Suppose that $H^\phi = 0$, then $(4.3)$ implies
\begin{equation}
	g(X, X)\triangle f - H^f(X, X) = \{\tilde{\lambda}(n - 1) + r\}fg(X, X),
\end{equation}
for all smooth vector fields $X$ on $M$. Now assume that $g(X, X)\triangle f - H^f(X, X) \geq 0$, then $(4.4)$ provides
\begin{equation}
	\tilde{\lambda} \leq -\frac{r}{n - 1}
\end{equation}
Thus, considering $(4.5)$ together with $proposition(2.1)$, we get the desired result.
\end{proof}

\n Now, considering $theorem(4.2)$ with $proposition(2.2)$ we can state the following theorem.
\begin{theorem}
	Let $\tilde{M} =  M \times_f \Re$ be a static perfect fluid spacetime with $(M, g)$ a Ricci flat Riemannian manifold and dim$M \geq 2$, and $(\tilde{g}, \tilde{\lambda}, \tilde{\phi})$ is an almost gradient Ricci soliton on $\tilde{M}$ with $H^\phi = 0$, where $\phi = \tilde{\phi}\vert_M$ then $\tilde{M}$ satisfies the strong energy condition if 
	\begin{equation}\nonumber
		\tilde{\lambda} \leq \frac{r}{n - 1}
	\end{equation}
\end{theorem}

\section{Almost Gradient Ricci Solitons on Base Manifold of Static Perfect fluid spacetime}
In this section, we study an almost gradient Ricci soliton structure on base manifold $M$ of static perfect fluid spacetime $\tilde{M} =  M \times_f \Re$ with potential function as warping function and soliton constant $\lambda$, defined as
\begin{equation}
	H^f(X_1, X_2) + S(X_1, X_2) + \lambda g(X_1, X_2) = 0,
\end{equation}
for all smooth vector fields $X_1, X_2$ on $M$, where $H^f$ denotes the Hessian of smooth function $f$ on $M$, $S$ is the Ricci tensor of $g$ and $\lambda$ is a smooth function on $M$.
	\begin{theorem}
		If the base manifold $M$ of a SPF spacetime $\tilde{M} =  M \times_f \Re$ is admitting an almost gradient Ricci soliton then it becomes an Einstein manifold.
	\end{theorem}
\begin{proof}
	Let the base manifold $(M^n, g)$ of a SPF spacetime $(\tilde{M}^{n + 1}, \tilde{g})$ admits an almost gradient Ricci soliton. Using $(2.3)$ in $(5.1)$, we get
	\begin{equation}
		f(S - \frac{r}{n}g)(X_1, X_2) + \frac{\triangle f}{n}g(X_1, X_2) + S(X_1, X_2) + \lambda g(X_1, X_2) = 0,
	\end{equation}
for all smooth vector fields $X_1, X_2$ on $M$. Now using $(2.4)$ and $(5.2)$ together, we obtain
\begin{equation}
	S(X_1, X_2) = \frac{1}{(f + 1)}\left\{\frac{r - 2\rho}{2(n - 1)}f - \lambda\right\}g(X_1, X_2),
\end{equation}
for all smooth vector fields $X_1, X_2$ on $M$. This completes the proof.
\end{proof}
\begin{remark}
	\item[(i)] Tracing $(5.3)$ over $X_1$ and $X_2$, the soliton constant of almost gradient Ricci soliton is 
	\begin{equation}
		\lambda = \frac{(r - 2\rho)}{2(n - 1)}f - \frac{r(f + 1)}{n}.
	\end{equation}
\item[(ii)] Since for Stiff matter fluid, $\sigma = \rho$ and for spatial factor $\sigma = \frac{r}{2}$, so the isotropic pressure $\rho = \frac{r}{2}$. Also, using $\rho = \frac{r}{2}$ in $(5.4)$, we obtain
\begin{equation}
	\lambda = -\frac{r(f + 1)}{n}.
\end{equation}
Thus, we can write:
If the base manifold $(M^n, g)$ of a SPF spacetime $(\tilde{M}^{n + 1}, \tilde{g})$ admitting an almost gradient Ricci soliton is a stiff matter fluid then the soliton is shrinking, provided the scalar curvature is positive.
\end{remark}

	\begin{definition}$[9]$
		A perfect fluid spacetime is said to have the timelike convergence condition if the Ricci tensor $S$ of the spacetime satisfies the condition $S(X, X) > 0,$ for all timelike vector fields $X$ on $M$.
	\end{definition}
\begin{theorem}
	If the base manifold $(M^n, g)$ of a SPF spacetime $(\tilde{M}^{n + 1}, \tilde{g})$ is admitting an almost gradient Ricci soliton then it obeys timelike convergence condition if and only if the soliton constant $\lambda > \frac{2\rho - r}{2(n - 1)}f.$
\end{theorem}
\begin{proof}
	Now, taking $X_1 = X_2 = \xi$ in $(5.3)$, we get
	\begin{equation}
		S(\xi, \xi) = \frac{1}{(f + 1)}\left\{\frac{r - 2\rho}{2(n - 1)}f - \lambda\right\}g(\xi, \xi).
	\end{equation}
	Consider the base manifold $(M^n, g)$ of a SPF spacetime $(\tilde{M}^{n + 1}, \tilde{g})$ is obeying timelike convergence condition $(i.e. S(\xi, \xi) > 0)$, then $(5.6)$ provides
	\begin{equation}\nonumber
		\frac{1}{(f + 1)}\left\{\lambda + \frac{r - 2\rho}{2(n - 1)}f \right\} > 0.
	\end{equation}
	which implies
	\begin{equation}\nonumber
		\lambda > \frac{2\rho - r}{2(n - 1)}f.
	\end{equation}
This completes the proof.
\end{proof}
\begin{proposition}
	If the base manifold $(M^n, g)$ of a SPF spacetime $(\tilde{M}^{n + 1}, \tilde{g})$ admits an unit torseforming vector field then the curvature tensor and Ricci tensor satisfies the relations:
	\begin{equation}\nonumber
		R(X_1, X_2)\xi = \eta(X_2)X_1 - \eta(X_1)X_2,
	\end{equation}
and \begin{equation}\nonumber
	S(X_1, \xi) = (n - 1)\eta(X_1).
\end{equation}
\end{proposition}
\begin{proof}
	Since the base manifold $(M^n, g)$ of a SPF spacetime $(\tilde{M}^{n + 1}, \tilde{g})$ admits a unit torseforming vector field $\xi$, therefore from $(2.19)$, we infer
	\begin{equation}
		\nabla_{X_1}\xi = X_1 + \eta(X_1)\xi,
	\end{equation}
for all smooth vector field $X_1$ on $M$. Now, taking covariant derivative of $(5.7)$ along $X_2$, we get
\begin{equation}
	\nabla_{X_2}\nabla_{X_1}\xi = \nabla_{X_2}X_1 + (\nabla_{X_2}\eta)(X_1)\xi + \eta(\nabla_{X_2}X_1)\xi + \eta(X_1)\nabla_{X_2}\xi,
\end{equation}
for all smooth vector field $X_1, X_2$ on $M$. Using $(5.8)$ in the definition of Riemannian curvature tensor
\begin{equation}\nonumber
	R(X_1, X_2)\xi = \nabla_{X_1}\nabla_{X_2}\xi - \nabla_{X_2}\nabla_{X_1}\xi - \nabla_{[X_1, X_2]}\xi,
\end{equation}
we obtain
\begin{equation}
		R(X_1, X_2)\xi = \eta(X_2)X_1 - \eta(X_1)X_2,
\end{equation}
for all smooth vector field $X_1, X_2$ on $M$. Further contracting $X_2$ from $(5.9)$, we get
\begin{equation}
	S(X_1, \xi) = (n - 1)\eta(X_1),
\end{equation}
for all smooth vector field $X_1$ on $M$. This finishes the proof.
\end{proof}
\begin{theorem}
	If the base manifold $(M^n, g)$ of a SPF spacetime $(\tilde{M}^{n + 1}, \tilde{g})$ is admitting an almost gradient Ricci soliton and the basic vector field $\xi$ is torseforming vector field then either the torseforming vector field is an orthogonal vector field, or the soliton constant is,
	\begin{equation}
		\lambda = \frac{(r - 2\rho)f - 2(n - 1)^2(f + 1)}{2(n - 1)}.
	\end{equation}
\end{theorem}
\begin{proof}
	Setting $X_2 = \xi$ in $(5.3)$ entails that
	\begin{equation}
		S(X_1, \xi) = \frac{1}{(f + 1)}\left\{\frac{r - 2\rho}{2(n - 1)}f - \lambda\right\}\eta(X_1),
	\end{equation}
for all smooth vector field $X_1$ on $M$. Now, $(5.10)$ and $(5.12)$ together provides
\begin{equation}
	\left[(n - 1) + \frac{1}{(f + 1)}\left\{\frac{2\rho - r}{2(n - 1)}f + \lambda\right\}\right]\eta(X_1) = 0,
\end{equation}
for all smooth vector field $X_1$ on $M$. This gives either $\eta(X_1) = 0$, or the soliton constant is given by $(5.11)$. This completes the proof.
\end{proof}
\begin{remark}
	Since for Stiff matter fluid, $\sigma = \rho$ and for SPF spacetime $\sigma = \frac{r}{2}$, so the isotropic pressure $\rho = \frac{r}{2}$. Also, using $\rho = \frac{r}{2}$ in $(5.11)$, we obtain
	\begin{equation}
		\lambda = -(n - 1)(f + 1).
	\end{equation}
Thus, we can write:
If the base manifold $(M^n, g)$ of a SPF spacetime $(\tilde{M}^{n + 1}, \tilde{g})$ admitting an almost gradient Ricci soliton is a stiff matter fluid and the basic vector field is torseforming vector field then the soliton is always shrinking.
\end{remark}
\begin{theorem}
	Let the base manifold $(M^n, g)$ of a SPF spacetime $(\tilde{M}^{n + 1}, \tilde{g})$ admits a steady gradient Ricci soliton. If the base manifold is Ricci symmetric with constant scalar curvature then either warping function is constant, or it becomes the stiff matter fluid.
\end{theorem}
\begin{proof}
	Let the base manifold $(M^n, g)$ of a SPF spacetime $(\tilde{M}^{n + 1}, \tilde{g})$ is Ricci symmetric then
	\begin{equation}
		(\nabla_{X_1}S)(X_2, X_3) = 0,
	\end{equation}
for all smooth vector fields $X_1, X_2, X_3$ on $M$. Now, taking covariant derivative of $(4.3)$, we obtain
\begin{equation}
	(\nabla_{X_1}S)(X_2, X_3) = \left\{\frac{r - 2\rho}{2(n - 1)} - \lambda\right\}\frac{(X_1f)}{(f + 1)^2}g(X_2, X_3),
\end{equation}
for all smooth vector fields $X_1, X_2, X_3$ on $M$. Comparing $(5.15)$ and $(5.16)$, we infer
\begin{equation}
	\left\{\frac{r - 2\rho}{2(n - 1)} - \lambda\right\}\frac{(X_1 f)}{(f + 1)^2} = 0,
\end{equation}
for all smooth vector field $X_1$ on $M$. This provides either $(X_1 f) = 0$, or the soliton constant is,
\begin{equation}
	\lambda = \frac{r - 2\rho}{2(n - 1)}.
\end{equation}
Since the soliton is steady gradient Ricci soliton i.e. $\lambda = 0$, so $(5.18)$ implies the base manifold $(M^n, g)$ of a SPF spacetime $(\tilde{M}^{n + 1}, \tilde{g})$ is stiff matter fluid. Hence, the proof is completed.
\end{proof}

\n Now, we consider base manifold $(M^n, g)$ of a SPF spacetime $(\tilde{M}^{n + 1}, \tilde{g})$ is weakly Ricci symmetric.
\begin{theorem}
	Let the base manifold $(M^n, g)$ of a SPF spacetime $(\tilde{M}^{n + 1}, \tilde{g})$ admits a gradient Ricci soliton and the basic vector field is torseforming vector field. If the base manifold is weakly Ricci symmetric with constant scalar curvature then
	\begin{equation}
		\xi(f) = \frac{2(n^2 - 1)(f + 1)r}{2\lambda(n - 1) - (r - 2\rho)f}.
	\end{equation}
\end{theorem}
\begin{proof}
	Let us assume that base manifold $(M^n, g)$ of a SPF spacetime $(\tilde{M}^{n + 1}, \tilde{g})$ is weakly Ricci symmetric and admits a gradient Ricci soliton. Then from $(5.1)$ we get
	\begin{equation}
		\nabla_{X_1}Df + QX_1 + \lambda X_1 = 0,
	\end{equation}
for all smooth vector field $X_1$ on $M$. Using $(5.20)$ in the definition of Riemannian curvature tensor
\begin{equation}\nonumber
	R(X_1, X_2)Df = \nabla_{X_1}\nabla_{X_2}Df - \nabla_{X_2}\nabla_{X_1}Df - \nabla_{[X_1, X_2]}Df,
\end{equation}
we obtain
\begin{equation}
	R(X_1, X_2)Df = (\nabla_{X_2}Q)X_1 - (\nabla_{X_1}Q)X_2,
\end{equation}
for all smooth vector fields $X_1, X_2$ on $M$. Since base manifold $(M^n, g)$ of a SPF spacetime $(\tilde{M}^{n + 1}, \tilde{g})$ is weakly Ricci symmetric, therefore $(2.17)$ provides
\begin{equation}
	QX_1 = -r\eta(X_1)\xi,
\end{equation}
for all smooth vector fields $X_1$ on $M$. Now, taking covariant derivative $(5.22)$ along $X_2$ and using $(2.19)$, we get
\begin{equation}
	(\nabla_{X_2}Q)X_1 = -r\{g(X_1, X_2)\xi + \eta(X_1)X_2 + 2\eta(X_1)\eta(X_2)\xi\}.
\end{equation}
Interchanging $X_1$ and $X_2$ in above equation, we obtain
\begin{equation}
	(\nabla_{X_1}Q)X_2 = -r\{g(X_1, X_2)\xi + \eta(X_2)X_1 + 2\eta(X_1)\eta(X_2)\xi\}.
\end{equation}
Taking $(5.21)$, $(5.23)$ and $(5.24)$ together, we acquire
\begin{equation}
	R(X_1, X_2)Df = r\{\eta(X_2)X_1 - \eta(X_1)X_2\},
\end{equation}
for all smooth vector fields $X_1, X_2$ on $M$. Now, contracting $X_1$ from $(5.25)$ and for given local orthonormal frame $\{e_1, e_2,...,e_n\}$, we get
\begin{equation}
	S(X_2, Df) = r\{n\eta(X_2) + g(X_2, e_i)\},
\end{equation}
for all smooth vector field $X_2$ on $M$. Putting $X_1 = Df$ in $(4.3)$, we obtain
\begin{equation}
	S(X_2, Df) = \frac{1}{(f + 1)}\left\{\frac{r - 2\rho}{2(n - 1)}f - \lambda\right\}g(X_2, Df),
\end{equation}
for all smooth vector field $X_2$ on $M$. Comparing $(5.26)$ and $(5.27)$ and setting $X_2 = \xi$, we get
\begin{equation}
	\frac{1}{(f + 1)}\left\{\frac{r - 2\rho}{2(n - 1)}f - \lambda\right\}g(\xi, Df) + r(n + 1) = 0,
\end{equation}
which yields $(5.19)$. This completes the proof.
\end{proof}
\begin{theorem}
	Let the base manifold $(M^n, g)$ of a SPF spacetime $(\tilde{M}^{n + 1}, \tilde{g})$ admits an almost gradient Ricci soliton. If the base manifold is weakly Ricci symmetric then the soliton constant is,
	\begin{equation}
		\lambda = \frac{r}{n}f - \frac{\triangle f}{n} + r(f + 1).
	\end{equation}
\end{theorem}
\begin{proof}
	Taking $(2.3)$ and $(5.1)$ together, we obtain
	\begin{equation}
		({\nabla}^2 f)(X_1, X_2) = \left(\frac{\triangle f - rf - n\lambda f}{n(f + 1)}\right)g(X_1, X_2),
	\end{equation}
for all smooth vector fields $X_1, X_2$ on $M$. Since the base manifold $(M^n, g)$ of a SPF spacetime $(\tilde{M}^{n + 1}, \tilde{g})$ is weakly Ricci symmetric, therefore using $(2.17)$ with $(5.1)$, we get
\begin{equation}
	({\nabla}^2 f)(X_1, X_2) = r\eta(X_1)\eta(X_2) - \lambda g(X_1, X_2),
\end{equation}
for all smooth vector fields $X_1, X_2$ on $M$. Comparing $(5.30)$ and $(5.31)$ and setting $X_2 = \xi$ gives $(5.29)$. This finishes the proof.
\end{proof}

\n Here, we are discussing almost gradient Ricci soliton on 4-dimensional half conformally flat base manifold $(M^n, g)$ of a SPF spacetime $(\tilde{M}^{n + 1}, \tilde{g})$. The following algebraic characterization of self-dual algebraic curvature tensors will be used, in the next theorem.
\begin{lemma}$[3]$
	Let $(V, <\cdot, \cdot>)$be an oriented four-dimensional inner product space of neutral signature.
	\begin{enumerate}
		\item[(i)] An algebraic curvature tensor $R$ is self-dual if and only if for any positively oriented orthonormal basis $\{e_1, e_2, e_3, e_4\}$
		\begin{equation}
			\mathcal{W}(e_1, e_i, x, y) = \omega_{ijk}\epsilon_j\epsilon_k\mathcal{W}(e_j, e_k, x, y),
		\end{equation}
		for any $x, y \in V$, for $i, j, k \in \{2,3,4\}$ and where $\omega_{ijk}$ represents the signature of the corresponding permutation.
		\item[(ii)] An algebraic curvature tensor $R$ is self-dual if and only if for a positively oriented pseudo-orthonormal basis $\{t,u,v,w\}$ (i.e. the non-zero inner products are$〈t,v〉 = 〈u,w〉 = 1$)and for every $x, y \in V$,
		\begin{equation}
			\mathcal{W}(t, v, x, y) = \mathcal{W}(u, w, x, y), \quad \mathcal{W}(t, w, x, y) = 0, \quad \mathcal{W}(u, v, x, y) = 0.
		\end{equation}
	\end{enumerate}
\end{lemma}
\begin{theorem}
	If 4-dimensional half conformally flat base manifold $(M^n, g)$ of a SPF spacetime $(\tilde{M}^{n + 1}, \tilde{g})$ admitting an almost gradient Ricci Soliton of neutral signature is stiff matter fluid then the soliton is expanding, steady, or shrinking according as the scalar curvature is positive, zero, or negative, provided $\Vert \nabla f \Vert = 0$.
\end{theorem}
\begin{proof}
	Since, ${\Vert \nabla f \Vert}^2 = g(\nabla f, \nabla f) = 0$, taking covariant derivative of it, we obtain
	\begin{equation}
		({\nabla}^2 f)(\nabla f, X_1) = 0,
	\end{equation} 
for all smooth vector field $X_1$ on $M$. Combining $(5.34)$ with $(5.1)$ and using definition of Ricci operator, we get
\begin{equation}
	Q(\nabla f) = \lambda \nabla f.
\end{equation}
Since the base manifold $(M^n, g)$ of a SPF spacetime $(\tilde{M}^{n + 1}, \tilde{g})$ is 4-dimensional half conformally flat space, so we can suppose that it is self dual. Then for pseudo-orthonormal frame $\{\nabla f, e_1, e_2, e_3\}$, equations $(5.33)$ and $(5.1)$ implies
\begin{equation}
	\mathcal{W}(\nabla f, e_2, X_1, X_2) = \mathcal{W}(e_1, e_3, X_1, X_2) = \mathcal{W}(\nabla f, e_3, X_1, X_2) = \mathcal{W}(e_1, e_2, X_1, X_2) = 0,
\end{equation}
for all smooth vector fields $X_1, X_2$ on $M$. Taking $X_2 = \nabla f$ in first equation of $(5.36)$, and using $(2.22)$ and $(5.35)$, we obtain
\begin{equation}
	\mathcal{W}(\nabla f, e_2, X_1, \nabla f) = \left(\frac{r - 4\lambda}{6} + \frac{\triangle f}{3f}\right)g(\nabla f, X_1),
\end{equation}
for all smooth vector field $X_1$ on $M$. Now, setting $X_2 = \nabla f$ in second equation of $(5.36)$, and using $(2.22)$ and $(5.35)$, we get
\begin{equation}
	\mathcal{W}(e_1, e_3, X_1, \nabla f) = 0,
\end{equation}
for all smooth vector field $X_1$ on $M$. Combining $(5.37)$ with $(5.38)$, we receive either $\lambda = \frac{r}{4} + \frac{\triangle f}{2f},$ or $g(\nabla f, X_1) = 0$ for all smooth vector field $X_1$ on $M$. But $g(\nabla f, X_1) = 0$ provides $\nabla f \neq 0$, which leads to contradiction. Thus, we obtain
\begin{equation}
	\lambda = \frac{r}{4} + \frac{\triangle f}{2f}.
\end{equation}
Now, taking $(2.4)$ in $(5.39)$, we get
\begin{equation}
	\lambda = \frac{5r - 8\rho}{12}.
\end{equation}
Since the base manifold is a stiff matter fluid, so putting $\rho = \frac{r}{2}$ in $(5.40)$, we obtain
\begin{equation}
	\lambda = \frac{r}{12}.
\end{equation}
Hence, the proof is completed.
\end{proof}
\begin{corollary}
	If anti-self dual isotropic four dimensional half conformally flat base manifold $(M^n, g)$ of a SPF spacetime $(\tilde{M}^{n + 1}, \tilde{g})$ admitting an almost gradient Ricci soliton is stiff matter fluid then then the soliton is steady gradient Ricci soliton.
\end{corollary}

	\n \textbf {Acknowledgement} : Tarun Saxena gratefully acknowledges the financial support provided by Department of Science and Technology (DST), New Delhi, India [INSPIRE-Fellowship Code No. IF200228].
	
	\n \textbf {Author Contribution} : Akhilesh Yadav and Tarun Saxena wrote the main manuscript text together.
	
	\n \textbf {Data Availability Statement} : Data sharing is not applicable to this article as no datasets were generated or analysed during the study.
	
	\n \textbf {Conflict of Interest} : The authors declare that they have no conflict of interest.

	 \newpage
	
	\noindent\author{Akhilesh Yadav}\\
	\date{Department of Mathematics, Institute of Science, \\Banaras Hindu University, Varanasi-221005, India}\\
	\maketitle {\noindent E-mail: akhilesha68@gmail.com}

	\noindent\author{Tarun Saxena}\\
	\date{Department of Mathematics, Institute of Science, \\Banaras Hindu University, Varanasi-221005, India}\\
	\maketitle {\noindent E-mail: tarunsaxena254@gmail.com}

\end{document}